\documentclass[11pt]{article}

\usepackage{graphicx}
\usepackage{amsthm}
\usepackage{amsmath}
\usepackage{amssymb}
\usepackage{epstopdf}

\theoremstyle{plain}
\newtheorem{theorem}{Theorem}
\newtheorem{lemma}{Lemma}
\newtheorem{proposition}[theorem]{Proposition}

\theoremstyle{definition}

\def\@seccntformat#1{Section\ \csname the#1\endcsname\quad}

\def\be{\begin{enumerate}}
\def\ee{\end{enumerate}}
\def\bi{\begin{itemize}}
\def\ei{\end{itemize}}

\title{A note on closed 3-braids \footnote{this article is an expanded version of the preprint  ``A note on transversal knots which are closed 3-braids", arXiv math/0703669.}}
\author{Joan S. Birman\footnote{partially supported by the US National Science Foundation Grant DMS 0405586.} \ and William W. Menasco\footnote{partially supported by the US National Science Foundation Grant DMS 0306062}}
\date{May 13, 2008}

\begin{document}
\maketitle

\centerline{Dedicated to the memory of Xiao-Song Lin}

\begin{abstract}  This is a review article about knots and links of braid index 3. Its goal is to gather together, in one place, some of the tools that are special to knots and links of braid index 3, in a form that could be useful for those who have a need to calculate, and need to know precisely all the exceptional cases.   
\end{abstract}

Knots and links which are closed 3-braids are a very special class. Like 2-bridge knots and links, they are simple enough to admit a complete classification. At the same time they are rich enough to serve as a source of examples on which, with luck, a researcher may be able to test various conjectures. As an example, non-invertibility is a property of  knots which is fairly subtle, yet  we know (see Theorem~\ref{T:invertibility} below) precisely which links of braid index 3 are and are not invertible.  Therefore 3-braids could be very useful for testing potential invariants that might detect non-invertibility.  

The purpose of this note is to gather together, in one place, some of the tools that are special to knots and links of braid index 3, in a form that could be useful for those who have a need to calculate, and need to know precisely all the exceptional cases.    It includes the work in an earlier preprint \cite{BM-2006} which was restricted to 3-braid representatives of transversal knots.  The preprint \cite{BM-2006} was written when \cite{NOT} was in preparation, and the authors of \cite{NOT} asked us about low crossing examples to test their then-new invariants of transversal knots.   

%%%%%%%%%%%%%%%%%%%%%%%%%%%%%%%%
%%%%%%%%%%%%%%%%%%%%%%%%%%%%%%%%
%%%%%%%%%%%%%%%%%%%%%%%%%%%%%%%%
\section{Three basic theorems}\label{S:Three basic theorems}

In this section we present the basic lemmas and theorems about links that are closed three-braids. 
Initially, we will use the classical presentation for the braid group $B_3$:
\begin{equation}\label{E:classical presentation}
<\sigma_1,\sigma_2; \sigma_1\sigma_2\sigma_1 = \sigma_2\sigma_1\sigma_2>
\end{equation}
In the presentation (\ref{E:classical presentation}) the generator $\sigma_i$ is a positive interchange of strands $i$ and $i+1$, where $i=1$ or $2$. 

Before we can state our first theorem, we need a definition. Let $\mathcal L$ be a closed 3-braid, with  representative $L\in B_3$. Then:
\bi
\item  $\mathcal L$ is said to {\it admit a flype} if the conjugacy class  of $L$ has a representative which is in one of the special forms illustrated in each of the four sketches in Figure~1.  In particular its associated (open) braids are conjugate to:
\begin{equation} \label{E:flype-admissible form} 
L = \sigma_1^u\sigma_2^v\sigma_1^w\sigma_2^\epsilon,  \  \  {\rm for \ \ some} \ \ u,v,w,\epsilon \in \mathbb N,  \ \ \epsilon = \pm 1.
\end{equation}
\item  A flype  is  {\it positive} or {\it negative}, according as $\epsilon$ is $+1$ or $-1$.  
\item   A flype is {\it non-degenerate} when the braids $L= \sigma_1^u\sigma_2^v\sigma_1^w\sigma_2^\epsilon$  and $L' = \sigma_1^w\sigma_2^v\sigma_1^u\sigma_2^\epsilon$ are in distinct conjugacy classes.  We are interested only in non-degenerate flypes. 
\ei
\begin{figure}[htpb!]
\label{F:flype}
\centerline{\includegraphics[scale=0.8] {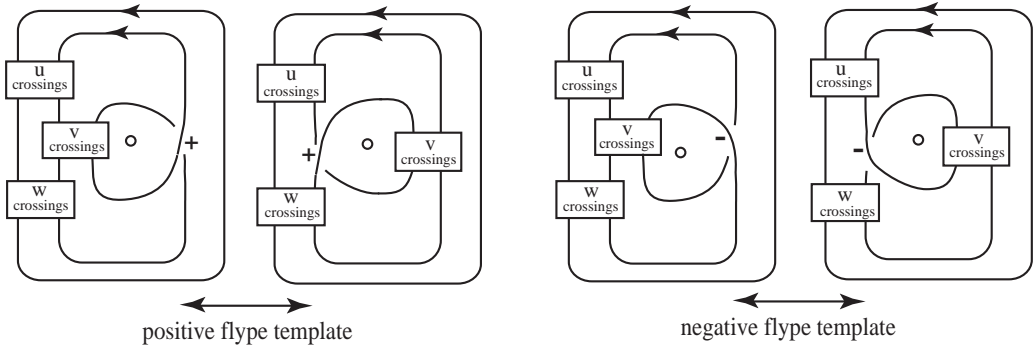}}
\caption{Positive and Negative 3-braid flypes}
\end{figure}

\begin{theorem} {\bf The classification theorem \cite{BM-SLVCB-III}:}
\label{T:topological}
Let $\mathcal{L}$ be a link type which is represented by the closure of a 3-braid $L$.  Then one of the following holds:
\be
\item  [{\rm (1)}]  $\mathcal L$ has braid index 3 and every 3-braid which represents $\mathcal L$ is conjugate to $L$.
\item [{\rm (2)}]  $\mathcal L$ has braid index 3, and is represented by exactly two distinct conjugacy classes of closed 3-braids. This happens if and only if the conjugacy class of  $L$ contains a braid which admits a non-degenerate flype.
\item [{\rm (3)}] The braid $L$ is conjugate to $\sigma_1^k\sigma_2^{\pm 1}$ for some $k\in\mathbb{N}$. These are precisely the links which are defined by closed 3-braids, but have braid index less than 3.
\ee
\end{theorem}
{\bf Remark:}  \  It is immediately clear that in order to use Theorem~\ref{T:topological}, the reader will need definitive tests for deciding when two 3-braids are conjugate, and when a conjugacy class admits a non-degenerate flype, and when a conjugacy class contains $\sigma_1^k\sigma_2^{\pm 1}$ for some $k$.   These matters will be discussed in $\S$\ref{S:conjugacy} below.

Let $\mathcal K$ be a knot. The standard knot tables classify \underline{unoriented} knots.  A sharper classification assumes that we fix an orientation on both $S^3$ and on a representative $\overrightarrow{K}$.  Let $\overleftarrow{K}$ denote $\overrightarrow{K}$ with its orientation reversed.  We then define $\mathcal K$ to be {\it invertible} if  there is an orientation-preserving homeomorphism $h$ of $S^3$ such that $h(\overrightarrow{K})$ is smoothly isotopic to $\overleftarrow{K}$.  Our second classification theorem treats the question of when an oriented link of braid index 3 fails to be invertible. As was mentioned earlier, this is an interesting issue, because at this writing there is no effective test that recognizes whether an arbitrary link is or is not invertible.  Using the standard 
presentation (\ref{E:classical presentation}) for $B_3$, we note that if $\overrightarrow{K}$ is defined by a cyclic word in the standard generators, then $\overleftarrow{K}$ will be defined by the same cyclic word, read backwards.  

\begin{theorem} {\bf The invertibility theorem \cite{BM-SLVCB-III}:} 
\label{T:invertibility}  Let $\mathcal K$ be a link of braid index 3 with oriented 3-braid representative $\overrightarrow{K}$. Then $\mathcal K$ is non-invertible if and only if $\overrightarrow{K}$ and $\overleftarrow{K}$ are in distinct conjugacy classes, and the class of $\overrightarrow{K}$ does not contain a representative which admits a non-degenerate flype.
\end{theorem}

{\bf Proof:}  If $\overrightarrow{K}$ admits a non-degenerate flype, then the isotopy of $S^3$ which realizes the flype takes $\overrightarrow{K} \to \overleftarrow{K}$, so that $\mathcal K$ is invertible.   If $\overrightarrow{K}$ does not admit a non-degenerate flype, then by Theorem~\ref{T:topological} the class of  $\overrightarrow{K}$ is the only class that represents the oriented link $\mathcal K$. The assertion follows.  $\|$

A knot is said to be {\it transversal} if it is everywhere transverse to the standard contact structure in $\mathbb R^3$.  A {\it transversal knot type} is the equivalence class of a transversal knot, under isotopies  which are restricted so that every intermediate representative is transversal.  Our third classification theorem relates to transversal knot types that are closed 3-braids. It uses Theorem~\ref{T:topological} and is related to Theorem~\ref{T:invertibility}, however in Theorem~\ref{T:transversal} below we are interested explicitly in the conjugacy classes which admit a non-degenerate negative flype and do not at the same time admit a positive flype, because a positive flype can be realized by a transversal isotopy.  Recall that a transversal knot type is {\it transversally simple} if it is determined completely by its topological link type and its Bennequin invariant.  

\begin{theorem}  {\bf The transversal classification theorem \cite{BM-stab-II}:}
\label{T:transversal}  
Let $\mathcal{TK}$ be a transversal knot type which is represented by a closed 3-braid 
$K_1 = \sigma_1^u\sigma_2^v\sigma_1^w\sigma_2^{-1}$.  Assume that $\mathcal K_1$ does not also admit a positive flype.  Then the braids
$K_1$ and 
$K_2 = \sigma_1^w\sigma_2^v\sigma_1^u\sigma_2^{-1}$  close to distinct transversal knot types $\mathcal{TK}_1, \mathcal{TK}_2$ that are topologically invertible but are not transversally invertible.  In particular, $\mathcal{TK}_1$ and $\mathcal {TK}_2$ are not transversally simple. 
\end{theorem}

%%%%%%%%%%%%%%%%%%%%%%%%%%%%%%%%
%%%%%%%%%%%%%%%%%%%%%%%%%%%%%%%%
%%%%%%%%%%%%%%%%%%%%%%%%%%%%%%%%
\section{Solving the conjugacy problem in $B_3$}  \label{S:conjugacy} 
As noted above, to use Theorem~\ref{T:topological} the reader will need definitive tests for deciding (i) when two 3-braids are conjugate; (ii) when a conjugacy class is one of the exceptional ones which admits a non-degenerate flype; and (iii) when a conjugacy class contains $\sigma_1^k\sigma_2^{\pm 1}$ for some $k$.    In this section we discuss (i). 

The conjugacy problem in $B_3$ was solved long ago by O. Schreier \cite{Sc}, making use of the fact that $B_3$ is a free product with amalgamations of cyclic groups of orders 2 and 3. His work did not generalize to $B_n, \ n>3$, where the braid groups have a more complicated structure than they do when $n=3$.  

The first solution to the conjugacy decision and conjugacy search problems for all braid groups $B_n, \ n\geq 3,$ was discovered by F. Garside in 1965 \cite{Ga}.  Garside's solution was exponential in word length and in braid index, and that fact  has lead to major efforts to find a polynomial-time solution to the conjugacy decision and search problems in $B_n$, which continue at this writing.    
In 1992, motivated by the work of Garside \cite{Ga} and the later work of Birman and Menasco on 3-braids in \cite{BM-SLVCB-III}, P.J. Xu discovered in \cite{Xu} a new and fast solution to the conjugacy problem in the special case of $B_3$. Her work lead naturally to an algorithm to solve the shortest word problem in $B_3$, and so (by special properties of links that are closed 3-braids, established in \cite{BM-SLVCB-III}) to compute the genus of a knot of braid index 3.  One surprising feature of it, that was generalized later to all the groups $B_n, n\geq 2,$ in \cite{BKL} and thence to many Garside groups, was that it gave a `dual' Garside structure on the braid groups.  This fact  played an important role in the isolation of {\it Garside groups} as an interesting class \cite{DP, Deh}.   In describing the basic theorems about 3-braids we choose Xu's work over Schreier's because it generalizes to all $n$ and so may give insights that might otherwise be missed.  We also choose the band presentation (described below) and its Garside structure over the classical presentation (\ref{E:classical presentation}) and its Garside structure because it seems simpler, from a computational point of view.  None of the solutions to the conjugacy problem in $B_3$ are particularly suited to problems (ii) and (iii) above.   

\subsection{The band presentation}
The  {\it band presentation} of $B_3$ uses a cyclically ordered set of generators which includes the generators of (\ref{E:classical presentation}):
\begin{equation} \label{E:BKL presentation}
<a_1, a_2, a_3 \ | \  a_2a_1 = a_3a_2 = a_1a_3> , \  {\rm where}  \  a_1=\sigma_1, \ \ a_2 = \sigma_2, \ \ a_3 = \sigma_1^{-1}\sigma_2\sigma_1.
\end{equation}  
In terms of the new generators, a 3-braid $L$ admits a flype (see (\ref{E:flype-admissible form})) if its conjugacy class contains $a_1^ua_2^va_1^wa_2^{\pm 1}$ for some $u,v,w,\epsilon$.   The constraints on $u,v,w,\epsilon$ are as before.
In working with the presentation (\ref{E:BKL presentation}) subscripts are always defined mod 3.  Let $\delta = a_2a_1 = a_3a_2 = a_1a_3$. The following formulas follow immediately from the defining relations in (\ref{E:BKL presentation}).
\begin{equation} \label{E:working with delta} 
(i) \ \  a_i^{-k} = \delta^{-k} a_{i-k+2}a_{i-k+1}\cdots a_ia_{i+1} \ \ {\rm for \ all} \ \ k\geq  1, \ \ (ii) \ \ a_j \delta^m = \delta^m a_{j+m} \  \ {\rm for \ all} \ \ m\in \mathbb N.
\end{equation}
They will be used repeatedly in what follows.  Note that $\delta^3$ generates the center of $B_3$, because subscripts are defined mod 3.  

\subsection{The word problem:}
Xu's solution to the word problem builds directly on the analogous ideas in the work of Garside \cite{Ga}.  Choose any word $W$ in the symbols  $a_1,a_2,a_3$ and their inverses, representing an element $\beta\in B_3$.  Using (\ref{E:BKL presentation})(i), all inverses of generators can be traded for powers of $\delta^{-1}$.  Using (\ref{E:BKL presentation})(ii), all powers of $\delta$ may be collected, say, at the left, so that  one may rewrite $W$ in the form $\delta^kP$, where $k\in \mathbb N$ and $P$ is a {\it positive} word.  
Positive words contain powers of $\delta$ if and only if two adjacent letters $a_{i+1}a_i$ appear, so after a simple search we may change $\delta^kP$ to a word of the form 
$\delta^qa_{r_1}a_{r_2}\cdots a_{r_s}$ where $q$ is maximal and the positive word $a_{r_1}a_{r_2}\cdots a_{r_s}$ contains no adjacent letters of the form $a_{i+1}a_i$.  Rewriting our word in syllable form $X = \delta^qa_{\tau_1}^{s_1}\dots a_{\tau_t}^{s_t}$, we see that the subscripts must be strictly increasing.   The integer $q$ is the {\it power}, and the positive word  is the {\it tail}.  This solves the word problem in $B_3$.  

We give an example. The following set of moves bring the flype-admissible braid $W=a_1^{-2}a_2^{-3}a_1^5 a_2$, which is not in normal form,  to its normal form $\delta^{-4}a_2^2a_3a_1^5a_2$:  
Use (i) to modify $W   = a_1^{-2}a_2^{-3}a_1^5 a_2$  to  
 $(\delta^{-1}a_2)^2(\delta^{-1}a_3)^3a_1^5a_2$.  Now use (ii) to modify $(\delta^{-1}a_2)^2(\delta^{-1}a_3)^3a_1^5a_2$ to $\delta^{-2} a_1a_2  \delta^{-3 }a_1 a_2a_3 a_1^5a_2$ and then to $ \delta^{-5}a_1a_2 a_1a_2a_3a_1^5a_2^4$. Finally, observe that $\delta=a_2a_1$ is a subword of the tail of $ \delta^{-5}a_1a_2 a_1a_2a_3a_1^5a_2^4$, so we can move it to the left to obtain $\delta^{-4}a_2^2a_3a_1^5a_2.$  The subscripts in the tail are strictly increasing, so we are done.

 \subsection{The conjugacy problem for 3-braids}  \label{SS:the conjugacy problem}
The conjugacy problem is close-by.  The {\it summit set} of $W$, denoted SS($W)$,  is the subset of those elements in the conjugacy class of $W$ which, when represented by words in normal form, have maximal power.  It is, by definition,  invariant under conjugation by powers of $\delta$.  Even more, it must be finite, for the following reasons. The defining relations in the band presentation preserve the algebraic sum of the exponents of an element in $B_3$, and conjugation does too. Therefore  the algebraic sum  is a class invariant, and from this it follows that every time that the power is increased by conjugation, the tail goes down in letter length, until we reach the summit set.  As for the summit set, it must be finite because the letter length of the tail of a word in SS$(W)$ is finite and is an invariant. The summit set  is Xu's basic invariant of the conjugacy class of $W$.  This was the invariant in Garside's algorithm in \cite{Ga}, however in the special case of $B_3$, using the band presentation, Xu was able to determine SS$(W)$ rapidly, and choose a unique symbol $\Sigma^\star(W)$ with the property that two elements $W,W'$ are conjugate if and only if $\Sigma^\star(W) = \Sigma^\star(W')$.   Her unique symbol is defined as follows. Suppose that 
$X = \delta^m a_{i_1}^{l_1}a_{i_2}^{l_2}\cdots a_{i_t}^{l_t}$ is in the summit set of $W$.  She associates to $X$  the array $\Sigma(X) = (m,(l_1,l_2,\dots,l_t))$.  The integer $m$ is fixed (and is minimal for the conjugacy class.  Her representative $\Sigma^\star(W)$ is chosen so that $t$ is minimal, and so that among all sequences with the same $t$  it is lexicographically smallest. Thus, for example,  $(5,(2,2))<  (5,(1,2,1))$ and $(5,(1,3))<(5,(2,2))$.  She proves that $W, W'\in B_3$ are conjugate if and only if $\Sigma^\star(W) = \Sigma^\star(W')$. 
We illustrate the algorithm via an example.  

\underline{Step 1:}   The first part of the algorithm is to find a single element in SS$(W)$. Let us suppose that we are given an element $W = \delta^qa_{\tau_1}^{s_1}\dots a_{\tau_t}^{s_t} \in B_3$, which is in normal form but may or may  not be in its summit set. To test whether $W$ is in SS$(W)$,  we modify $W= \delta^qa_{\tau_1}^{s_1}\dots a_{\tau_t}^{s_t}$ by  `twisted cyclic permutation' of  the tail of $W$.  That is, we push the left-most syllable $a_{\tau_1}^{s_1}$ in the tail to the left, past $\delta^q$, using (\ref{E:working with delta}) and then move it to the right end of $W$.  This may be achieved by a special kind of conjugation that we call  a {\it tail move} on $W$: 
$$W  = \delta^qa_{\tau_1}^{s_1}\dots a_{\tau_t}^{s_t} \to (\delta^qa_{\tau_1}^{-s_1}\delta^{-q})(W)(\delta^qa_{\tau_1}^{s_1}\delta^{-q}) = \delta^qa_{\tau_2}^{s_2}\cdots a_{\tau_t}^{s_t}(\delta^qa_{\tau_1}^{s_1}\delta^{-q}) = X.$$
Xu proves that $W$ is in its summit set if and only if, after a single tail move, the power is not increased and the syllable length is not decreased. To understand why this is so, observe that the final syllable of $X$ is $\delta^qa_{\tau_1}^{s_1}\delta^{-q} = a_j^{s_1}$ for some $j=1,2,3$ and so there are 3 possibilities, namely $j = \tau_t +1, \tau_t$ or $\tau_t - 1$ (mod 3).  Xu proves that if $j=\tau_t+1$ or $\tau_t$, then $W$ and also $X$ are in normal form as written and are in SS($W)$. The first part of the algorithm ends.  For,  if $X$ is modified again by another tail move, then the modified $X$ will have the same power as $X$ and be in normal form, because  $W$ was in normal  form, so that $\tau_2 = \tau_1 +1$ (mod 3), and this property  is preserved by the second tail move. 

If $X$ is not in normal form, then the only possibility is that  $j=\tau_t - 1$,  in which case we obtain a new factor $\delta$ at the interface between the final letter in $W$ and the initial letter in the syllable that was permuted, because $\delta = a_{\tau_t}a_{\tau_t-1}$.  An example is given by $W = \delta^{-3}a_1^4a_2a_3^3a_1^2a_2$.  After a tail move it changes to $X = \delta^{-3}(a_2a_3^3a_1^2a_2)(a_1^4) = \delta^{-3}a_2a_3^3a_1^2\delta a_1^3 =  \delta^{-2}a_3a_1^3a_2^2a_1^3$, which has higher power  than $W$.    But in this example $X$ is not in normal form, because when the factor $\delta$ that appeared after the tail move  was moved to the left, an additional factor $\delta$ popped up at the interface between the final two syllables. Pushing the new factor $\delta$ to the left, we obtain 
$Y = \delta^{-2}a_3a_1^3a_2\delta a_1^2 = \delta^{-1} a_1a_2^3a_3a_1^2.$  The power has been increased again.  Normal form has been achieved, because the subscript array is increasing. 

We now test whether $Y$ in our example is in SS$(W)$, by applying a tail move. The answer is `yes'.  More generally, the process of repeated applications of tail moves, and passage to normal form, continues until we find an element whose power is unchanged after a tail move.  In our example that element is $Y$ and so $Y \in {\rm SS}(W)$. In a more general example, another cycling might lead to additional powers of $\delta$, but in the most general example the following holds: If a tail move does not change the power, then we have achieved summit power.  Note that the tail move may, nevertheless, increase syllable length. 

\underline{Step 2:} \  Having achieved summit power, Xu's method for choosing an invariant of the conjugacy class is the following: If $Y= \delta^ma_{\mu_1}^{l_1}\dots a_{\mu_k}^{l_k}$ belongs to SS$(W)$, one associates to $Y$ its {\it Xu symbol}  $\Sigma(Y) = (m; (l_1,l_2,\dots,l_k))$.  In our example, $\Sigma(Y) = (-1; (1,3,1,2))$.  

Note that every element in SS$(W)$ has the same summit power $m$, but elements need not have the same syllable length.  In particular, if in step 1 it had happened that, after a tail move on  $W$, we found ourselves in the case $j = \tau_t$,  then the syllable length would go down,  because in $X$ the permuted syllable $a_j^{s_1}$ could be amalgamated with  the final syllable $a_{\tau_t}^{s_t}$ of $W$.  Among all elements $Y \in {\rm SS}(W)$, her defining symbol $\Sigma^\star(W)$ is the one for which the sequence of exponents $(l_1,l_2,\dots,l_k)$ is lexicographically smallest.  In particular, the syllable length is minimized among elements in the summit set.  It is now easy to see that if a tail move leaves both the power and the syllable length invariant then we are almost done with step 2. 

We study $\Sigma(Y) = (-1; (1,3,1,2))$ in our example.  We seek the unique element in SS$(W)$ for which the tail symbol is lexicographically smallest. But now observe that the candidates are simply the 4 cyclic permutations of the tail symbol, i.e. (1,3,1,2), (3,1,2,1), (1,2,1,3). (2,1,3,1).  (In a more general case the tail symbol will have, say, u  entries, and so it will have $u$ cyclic permutations, all of which occur for some element in SS$(W)$.)   In our example $\Sigma^\star(Y) = (-1; (1,2,1,3))$ is the Xu invariant of the conjugacy class. We found it without calculating the entire summit set, in particular we ignored all the elements whose syllable length was bigger than 4, and we did not have t0 deal with the 3 possible subscript sequences of each element that has this Xu symbol. Therefore we are done.  $\|$

\section{Using Theorems 1,2, and 3}

{\bf Theorems 1 and 2}. \ \ We need two more results before we can make Theorems~\ref{T:topological} and \ref{T:invertibility}  into results that can be used computationally.  First, we need to be able to recognize when the closure of a 3-braid has braid index less than 3. Table I  gives the Xu symbols $\Sigma^\star(W)$  in the cases when the closed braid associated to $W$ has braid index less than 3.  Using it, we can decide when a closed 3-braid is a link of braid index 3, by computing its Xu-invariant and verifying whether it is or is not the symbol for a  3-braid of braid index less than 3.  In the table the symbol $1^{u}$ means $u$ syllables of length 1.  For example, Table I tells us that a 3-braid $X$ whose Xu-invariant is $(-5; (1^4, 2)) = (-5; (1, 1, 1, 1, 2))$ is conjugate to the 3-braid  $a_1^{-5}a_2$, which closes to the type (2,5) torus knot.  Therefore the closure of $X$ has braid index 2. 

\medskip

\centerline{
\begin{tabular}{ |c | c| }
\hline
$W$ & $\Sigma^\star(W)$  \\
\hline
\hline
$a_1^k a_2$,  \ $k = 1 \ {\rm and} \ 2$  &$(1; \ (k-1))$ \\ \hline
$a_1^k a_2$,  \ $k \geq 3$  &$(-1; \ (k+1))$ \\ \hline
 $a_1^k a_2^{-1}$, \   $k\geq 2$      &  $(-1; \ (1,k))$   \\ \hline
$a_1^k a_2$,  \ \ \ $k\leq -2$ & $(k; (1^{-k-1},2))$ \\  \hline
\ \ \ $a_1^k a_2^{-1}$, \   $k \leq -1$, \   & $(k; \ (1^{-k-1}))$    \\ \hline
$a_1^{-1}a_2$ and $a_1a_2^{-1}$  &  $(-1; (2))$    \\ \hline
$a_2$  &  $(0; (1))$    \\ \hline
$a_2^{-1}$  &  $(-1; (1))$ \\ \hline
 \end{tabular} }
\medskip
\centerline{\bf Table I: Conjugacy classes of links of braid index less than 3}

Next, if we are given two 3-braids $W,W'$ and find that they have distinct Xu-invariants, we cannot be sure that they determine distinct links until we determine whether they are flype-admissible.  
To answer this question, we need the following Lemma, which is due to K.H.Ko and S.J.Lee:
\begin{lemma}{\rm K.H.Ko and S.J.Lee, \cite{KL}}
\label{L:properties of flypes}  \ \ \  {\rm (i)} The flype  in {\rm (\ref{E:flype-admissible form})} is non-degenerate if and only if  $|v|\geq 2$, also $u,v+\epsilon,w$ are pairwise distinct, also  $u,w \not= 0, \epsilon,  2\epsilon$. 
{\rm (ii)} The conjugacy class of a negative flype-admissible braid $L$ contains a braid that admits a positive flype  if and only if  $u=1$ or $w = 1$, or $v = 2$.  
\end{lemma}
In Table II we give the Xu symbols for all non-degenerate flype-admissible braids, up to cyclic permutation of the tail symbol.  There are 16 cases, because the exponents $u,v,w,\epsilon$ can either be positive or negative. To keep track of the signs, replace $u,v,w$ by $\pm p, \pm q, \pm r$, where $p,q,r>0$.  Now observe that only  8 of the 16 cases are needed, because a given braid $K$ admits a non-degenerate flype if and only if the Xu-symbol for $K$ or for $K^{-1}$ coincides with one of the 8 symbols that are given in Table II. 

\centerline{
\begin{tabular}{|c | c |c |  }
\hline
$\epsilon$ & $W$ &$\Sigma^\star(W)$, up to cyclic permutation of tail symbol  \\
\hline
\hline
$pos$ & $a_1^pa_2^qa_1^ra_2$ &$(3 \ ; \  (p-2, \ q-1, \ r-2))$ \\ \hline
$$ & $a_1^{-p}a_2^qa_1^ra_2$  & $(1-p \ ; \  (q, \ r-1, 1^p))$    \\ \hline
$$ & $a_1^pa_2^{-q}a_1^ra_2$   &   $(2-q \ ; \ (p-1, \ 1^{q-1}, \  r-1))$  \\ \hline
$$ & $a_1^pa_2^qa_1^{-r}a_2$ & $(1-r \ ; \  (p -1, \ q, \ 1^r ))$  \\ \hline
$neg$ & $a_1^pa_2^qa_1^ra_2^{-1}$ &   $(0 \ ; \  (p, \ q-1, \ r))$   \\ \hline
$$ & $a_1^{-p}a_2^qa_1^ra_2^{-1}$ &    $(-p \ ; \ (q, \ r+1, \ 1^{p-2}))$    \\ \hline
$$ &$a_1^pa_2^{-q}a_1^ra_2^{-1}$  &     $(-q-1 \ ; \ (r+1, \ p+1, \  1^{q-1}))$       \\ \hline
$$ & $a_1^pa_2^{q}a_1^{-r}a_2^{-1}$ &       $(-r \ ; \ (p+1, \ q, \ 1^{r-2}))$  \\ \hline
\end{tabular} }
\medskip
\centerline{{\bf Table II: Conjugacy classes admitting non-degenerate flypes}}
\medskip
We remark that the data in Table II is  only given up to cyclic permutation of the tail symbol, because the lexicographically smallest choice  depends upon the relative numerical values of $p,q,r$. 

We now have essentially all the data needed to implement Theorems 1,2 and 3, however an additional issue arises if one wishes to use any of our three theorems on low crossing examples which admit non-degenerate flypes.  For example, one might wish to test a potential invariant which was difficult to compute on knots which had large crossing number.   We give, in Table III, a list of all of the knots of braid index 3 which admit non-degenerate negative flypes, and do not also admit positive flypes, and have braid crossing number $c_b = |u|+|v|+|w| + 1\leq 12$.  Note that $c_b$ is not the same as minimum crossing number for all possible projections of the knot associated to our closed braid (although for all but one of the examples in Table 1 the two are actually the same), rather it is crossing number for closed 3-braid representatives which are in flype position. 

The knots are identified, in Table III, using both classical notation (e.g as in \cite{Rolfsen}) and the newer notation of Dowker-Thistlethwaite (used in all the newer tables).  The former only exist through 10 crossings.

\

\centerline{
\begin{tabular}{|c|c | c | c | c |c |  }
\hline
$\mathcal K$(classical) & $\mathcal K$(Dowker-Thistlethwaite) & $\beta$ & $c_b$ &$(u,v,w)$ & $(w,v,u)$ \\
\hline
\hline
$8_{10}$ & $8_{a3}$ & -1 &  8 & (3,-2,2) & (2,-2,3) \\ \hline
$10_{47}$ & $10_{a15}$ &  +1  & 10 &(5,-2,2) & (2,-2,5) \\\hline
$10_{62}$ & $10_{a 41}$ & +1 &  10  & (3,-2,4) & (4,-2,3)\\\hline
$10_{48}$ & $10_{a79}$ & -3 &   10 & (3,-4,2) & (2,-4,3)\\ \hline
$10_{127}$ & $10_{n16}$ &  -9       &  10    & (-5,-2,2) & (2,-2,-5)\\\hline
$10_{143}$ & $10_{n  26}$ & -7 &  10  & (3,-2,-4) & (-4,-2,3) \\\hline
$***$ & $11_{a240}$ & +7  &   12 & (5,3,3)& (3.3.5)  \\ \hline
$***$ & $12_{a 146}$ & +3 &   12 & (7,-2,2)& (2,-2,7)  \\\hline
$***$ & $12_{a 369}$ & +3 &   12  & (5,-2,4) &  (4,-2,5) \\\hline
$***$ & $12_{a 576}$ & +3 &  12  & (3,-2,6) & (6,-2,3)  \\\hline
$***$ & $12_{a  835}$ &-5 &   12  & (3,-6,2)& (2,-6,3) \\\hline
$***$ & $12_{a  878}$ & -1 &    12  & (5,-4,2) & (2,-4,5)  \\\hline
$***$ & $12_{a  1027}$ & -1 &    12  & (3,-4,4) & (4,-4,3)  \\\hline
$***$ & $12_{a  1233}$ & +1 &   12   &(5,-3,3) & (3,-3,5)  \\\hline
$***$ & $12_{n 234}$ & -11 &    12 & (-7,-2,2)&  (2,-2,-7) \\\hline
$***$ & $12_{n  466}$ & -9 &   12  & (-5,-3,3) & (3,-3,-5)  \\\hline
$***$ & $12_{n  467}$ & -7 &   12  & (-3,-4,4) & (4,-4,-3)  \\\hline
$***$ & $12_{n  468}$ & -5 &    12  & (-3,-3,5) & (5,-3,-3) \\\hline
$***$ & $12_{n  472}$ & -15 &    12  & (-3,-4,-4)& (-4,-4,-3)  \\\hline
$***$ & $12_{n  570}$ & -9 &   12  & (-3,-5,3)& (3,-5,-3)    \\
\hline
\end{tabular} }

\medskip

\centerline{{\bf Table III:  A list of all examples of non-transversally simple }}

\centerline{{\bf  closed 3-braids with braid crossing number $c_b\leq 12$}}

To assemble the data in Table III we began by listing all possible numbers $(\pm|u|,\pm|v|,\pm|w|)$   with $c_b=|u|+|v|+|w|+1$, deleting those that do not satisfy the constraints of the two lemmas.  Next we used Knotscape \cite{Knot} to sort the survivors  into classes which had the same  topological knot type.  For each topological knot type, some of the examples had non-minimal braid crossing number, so we deleted those. For example, $(-3,4,2)$ appeared as an example with $c_b=10$, but  its conjugacy class also has an 8-crossing example, $(3,-2,2)$, so we deleted $(-3,4,2)$.  After that, for each topological knot type we divided the representing closed braids into conjugacy classes belonging to the braids before and after the flype.  In every case it turned out that there were two distinct  conjugacy classes in $B_3$, as must be the case, by Theorem~\ref{T:topological}.  Finally, we chose exactly one representative from each conjugacy class, whenever it happened that the list of survivors had more than one representative.   In this regard we note that, for example, $\sigma_1^{-5}\sigma_2^{3}\sigma_1^{-3}\sigma_2^{-1}\approx 
\sigma_1^{2}\sigma_2^{-2}\sigma_1^{-5}\sigma_2^{-1} \approx 
\sigma_1^{-3}\sigma_2^{-4}\sigma_1^{2}\sigma_2^{-1}$.  This shows that conjugate braids which admit flypes need not have the same braid crossing number $c_b$. In this case the first has braid crossing number  12 and the other two have crossing number $10$. There are examples where  all three have the same crossing number, and examples where one has minimal crossing number and the other two have higher crossing number.  This issue arises when we enumerate the low-crossing number examples.

The examples in Table III are given by specifying,  for each pair of transversal knot types, the topological knot type $\mathcal{K}$ in the standard tables, the  Bennequin number $\beta = u+v+w-4$, the minimum closed braid crossing number $c_b=|u|+|v|+|w|+1$, and the two closed 3-braid representatives ${\mathcal {TK}}_1, {\mathcal {TK}}_2$. The closed braids are determined by the triplets $(u,v,w)$ and $ (w,v,u)$, i.e. as negative flype-related braids $\sigma_1^u\sigma_2^v\sigma_1^w\sigma_2^{-1}$ and $\sigma_1^w\sigma_2^v\sigma_1^u\sigma_2^{-1}$.

\section{Some open questions and a survey of known results about links of braid index 3} \label{S:known and unknown}  In this section we survey the literature and discuss some unknown properties of knots and links of braid index 3.  Our choices are guided by what we know best, and by our wish to illustrate ways in which Theorems 1,2 and 3 have been and might be useful.  

\begin{enumerate}
\item  We made inquiries, but were unable to find a computer program that implements Xu's solution to the conjugacy problem in $B_3$. If anyone is interested, it would be useful to have such a program, if it is freely available to researchers. 

\item    In our judgment, the most striking feature of closed 3-braids that has not been used in any substantial way in the literature is that we know precisely which links of braid index 3 are invertible, and which are not.   In the early history of knot theory it was not known whether non-invertible knots exist. The first proof that they do exist is due to H. Trotter \cite{Tr}, who showed that if $|p|,|q|,|r|\geq 2$ are distinct, then the assumption that the pretzel knots $K(p,q,r)$ are invertible leads to a contradiction about the structure of their groups.  His indirect proof characterizes all the work that we know which treats non-invertibility. We do not know of any computable numerical invariant that detects non-invertible knots. Therefore, if anyone has a candidate for an invariant that detetcts non-invertibility, three-braids are an excellent source of examples, because `most' of them are non-invertible.  Note that all knots that admit a flype, degenerate or not, are invertible. Theorem 2 and Tables I and II tell us precisely how to find all non-invertible 3-braid examples.

\item A knot $\overrightarrow{K}$ is {\it strongly invertible} if there is a representative $\overrightarrow{K}\subset S^3$  and an involution $\iota:(S^3,\overrightarrow{K}) \to (S^3,\overleftarrow{K})$.   In the course of writing this note we observed the following:

\begin{proposition}
In the special case of 3-braids, if a knot type $\overrightarrow{\mathcal K}$ is invertible, then it is strongly invertible. \end{proposition}
Question: Is it possible that this is the case for arbitrary knots?   

\item  We discuss next a problem which turns out to be closely related to problem 2 above, although initially it does not seem to be. A knot is {\it transversal} if it is everywhere transverse to the standard contact structure in $\mathbb R^3$.  The {\it transversal knot type} of a transversal knot is its equivalence class under isotopies through transversal knots. Known invariants of transversal knot types are their topological knot type  and their Bennequin number and (a recent development at the time of this writing) certain invariants based upon Knot Floer Homology  and Khovanov Homology \cite{OST, NOT, Pl1,Pl2, Wu,HKP}.  The first examples of pairs of knots which have the same topological knot type and the same Bennequin invariant, but are transversally distinct,  were discovered via 3-braids \cite{BM-stab-II}.  The proof in \cite{BM-stab-II} is indirect, and at this writing nobody has found an invariant which distinguishes even one of those pairs, however caution is needed: Each pair is an example of  a pair of knots which are topologically invertible but not transversally invertible.  In view of Problem 1 above they are bound to be extremely subtle to detect.

\item  We discuss a different set of questions about 3-braids which has received much attention and which, like Problem 2 above, seems very difficult. When the Jones and HOMFLY polynomials were discovered \cite{Jones} it was not known whether they were or were not injective on topological link types, but examples taken from 3-braids  \cite{Bir-Jones} quickly showed that the answer is `no'.   However, at this writing we know very little about how good a job the Jones polynomial does.   It is therefore a natural question to ask whether one can arrive at some understanding of the kernel of the Jones polynomial in the special case of braid index 3.   This is especially tempting because when the braid index is 3 the Jones polynomial depends on just one small calculation: the trace of a certain $2 \times 2$ matrix.  Nevertheless, the prohibitively complicated formulas in \cite{Tak}  tell us that this is simply not a good question, and that simplifications are needed before we can hope for a start in understanding the kernel of the Jones polynomial of links that are closed 3-braids.    

An example is enlightening, even though it does not have much to do with 3-braids. In the class of torus links, two integers $r,s$ are all that is needed to distinguish one torus link from another.  On the other hand,  if one wishes to distinguish the Jones polynomial of one torus link  from that of another,  the matter is considerably more complicated. One has the following very complicated formula \cite{Jones2}: 
\begin{equation} \label{E:JP of torus link}
V_{r,s}(t) = \frac{t^{\frac{1}{2}(r-1)(s-1)}}{1-t^2}\sum _{l=0}^{n}\left( n \over l \right)t^{\frac{r}{n}(n-l)(1+\frac{s}{n}l)}(t^{(n-l)\frac{s}{n}} - t^{1 + l\frac{s}{n}}), \ \  {\rm where} \ \   n = gcd(r,s).
\end{equation}
It isn't even obvious that it gives a polynomial (although of course it does).  The reader will notice that, knowing the general form of (\ref{E:JP of torus link}), $V_{r,s}(t)$  is completely determined by the two integers $r$ and $s$.
\item Concluding this note, we point the reader, in no particular order, to some of the published work that has come to our attention about links of braid index 3:
\begin{enumerate}
\item K. Murasugi's monograph \cite{Mur} was pioneering, in the sense that it was the first study that isolated knots and links of braid index 3 as a special class worthy of study.  It is filled with detailed information and the author's tireless, error-free and very useful calculations.  Murasugi  used Schreier's solution to the conjugacy problem \cite{Sc}, but carried things a step further by partitioning Schreier's conjugacy classes into groups which were later realized to coincide with pseudo-Anosov, finite order and reducible braids. In this sense it was prescient, because it forsaw the Nielsen-Thurston trichotemy that has played such a major role in understanding mapping class groups of surfaces,  in the special case of 3-braids. This seemed so striking that many people tried to carry generalize Murasugi's work to $B_n$, but alas the generalizations did not seem to work.  One may also find, in \cite{Mur}, an example  of the flypes which play such an important role in Theorems 1,2 and 3.   
\item  Two algorithmic solutions to the shortest word problem for 3-braids appeared in the literature essentially simultaneously. They were due to P.J.Xu \cite{Xu} and M. Berger \cite{Berger}, and they used the band and classical presentations respectively. A year later R.D.Keever \cite{Keev} gave yet another solution, duplicating both of the preceding. (His paper appeared in the same journal as Xu's paper, a year later, and neither Keever nor anyone else seems to have realized the duplication!).  Of these three, only Xu's paper does more, as it not only also gives the pleasant solution to the conjugacy problem that we discussed in $\S$\ref{SS:the conjugacy problem} but also proves that the shortest word problem implies an algorithm for computing link genus. 
Very recently  Usman Ali \cite{Us} posted a paper on the arXiv that chooses a unique representative in SS$(W)$ algorithmically when $W$ is a 3-braid, but using the classical presentation rather than the band presentation of Xu.   All these results are of interest because it has been proved by Paterson and Razborov \cite {PR} that the shortest word problem in $B_n$, using the classical presentation, is NP complete.
\item In a separate paper \cite{Xu2} Xu uses her earlier work to compute the growth functions for the positive braid semigroup.
Neither \cite{Berger} nor \cite{Keev} consider that problem.  
\item In a different direction S. Orevkov \cite{Orev} discovered an algorithm for deciding whether links of 
braid index 3 are quasi-positive, a question which is of interest in the study of the topology of plane real algebraic curves. 
\item In \cite{Ni} Y. Ni used knot Floer homology to determine which links of braid index 3 are fibered.  He introduced there the interesting concept of `almost fibered'.  The question of when 3-braids are fibered was also studied by A. Stoimenow in \cite{St1}.  
\item  M.T. Lozano and J. Przytycki \cite{LP} classified horizontal incompressible surfaces in the complement of a link of braid index 3.  Using a different set of techniques, E. Finkelstein obtained related results in \cite{Fink}. 
\item D. Erle has computed the signature of links of braid index 3 in \cite{Erle}.
\item A. Stoimenow gives a classification of alternating links of braid index 3 in \cite{St2}.

\item W. Menasco  and X. Zhang prove, in \cite{MZ}, that knots which are closed 3-braids have Property P.

\item In \cite{Kin} E. Kin uses 3-braids to study a problem in dynamics. 

\end{enumerate}
\end{enumerate}

{\bf Acknowledgment:} We thank Sang-Jin Lee, Lenny Ng, Peter Ozsvath and Dylan Thurston for useful comments on the precurser of this manuscript, the preprint \cite{BM-2006}.  We thank Usman Ali for pointing out errors in Table I.

 \end{document}